\documentclass[12pt]{article}
\usepackage[utf8]{inputenc}
\date{\today}
\author{Edward Eriksson}
\usepackage{amsthm}
\usepackage{tabularx}
\usepackage{ upgreek }
\usepackage[normalem]{ulem}
\usepackage{caption}
\usepackage{mathrsfs}
\usepackage{enumitem}
\usepackage[sorting=none]{biblatex}
\usepackage{amsmath}

\usepackage{authblk}  
\usepackage{algorithm}
\usepackage{algcompatible}
\usepackage{algpseudocode}
\usepackage[utf8]{inputenc}
\usepackage[UKenglish]{babel}
\usepackage[nottoc,numbib]{tocbibind}
\usepackage{mathtools,amssymb,amsthm}
\usepackage{adjustbox,graphicx}
\usepackage{mathtools}
\usepackage{hyperref}
\usepackage{lineno}
\usepackage{float}
\usepackage{tikz}
\usetikzlibrary{calc}
\usetikzlibrary{shapes.geometric, arrows}
\usepackage{csquotes}
\usepackage{amsthm}
\newtheorem{proposition}{Proposition}
\theoremstyle{definition}
\newtheorem{lemma}[proposition]{Lemma}

\newtheorem{theorem}[proposition]{Theorem}
\newtheorem{corollary}[proposition]{Corollary}
\newtheorem{definition}[proposition]{Definition}
\newtheorem{remark}[proposition]{Remark}
\newtheorem{example}[proposition]{Example}

\title{Edge Exchangeable Graphs:\\ Connectedness, Gaussianity and Completeness}
\author[1]{Edward Eriksson}
\affil[1]{Max Planck Institute for Mathematics in the Sciences, Leipzig}
\date{}

\begin{document}

\maketitle
\begin{abstract}
We characterize some asymptotic properties of edge exchangeable random graphs in terms of the measure used to generate them. In particular, we give a necessary and sufficient condition for eventual forever connectedness, a sufficient condition for asymptotic normality of the vertex count, and a necessary and sufficient condition for the produced graph to be eventually forever essentially complete. 
\end{abstract}
\section{Introduction} \label{Introduction}
\subsection{Model}

Following \cite{janson_2017_on,crane_2018_edge,cai_2016_edgeexchangeable,luo_2023_anomalous,li_2021_truncated} among others, we study models of the following type. We start with an empty (no vertices, no edges) graph and a probability measure $\mu$ on $\mathbb{N}_2$, the set of unordered pairs of distinct natural numbers. We think of the natural numbers as the set of possible vertices and $\mathbb{N}_2$ as the set of possible edges. Define a (multi-)graph-valued process by repeatedly sampling an edge according to $\mu$ and adding it, along with any needed vertices, to the graph. We also consider a Poissonized version where we sample from $\mu$ when there is an event in a Poisson process of unit intensity.

\begin{definition}
For $n \in \mathbb{N}$ let $e_n \sim \mu$, all independent. Define the multiset $E_n := \{e_1,\ldots,e_n\}$, $V_n := \cup_{i=1}^ne_i$ and
\begin{align*}
    G_n := (V_n,E_n).
\end{align*}
\end{definition}
For a set $A$ and a non-negative integer $k$ define $kA$ to be the multiset which consists of $k$ copies of each element of $A$.
\begin{definition}
For each $e \in \mathbb{N}_2$ let $N_e(t)$ be a Poisson process, all independent, with intensity $\mu_e := \mu(\{e\})$.  Let $E_t := \sum_{e \in \mathbb{N}_2}N_{e}(t)\{e\}$, where the sum is a multiset sum. Let  $V_t = \cup_{e \in E_t}e$ and finally set
$$G_t := (V_t,E_t).$$
\end{definition}

We will be interested in properties, such as connectedness and the number of vertices, which are invariant under identifying parallel edges, which we may therefore do. We will use the same notation for the graph-valued processes obtained after identifying parallel edges. Since the Poisson process has exponentially distributed waiting times we may equivalently think of the continuous time graphs as having edge $e$ added at time $\tau_e$ where the $\tau_e$ are independent exponential random variables with intensity $\mu_e$. By construction we have disallowed loops but this is largely a matter of convenience and elegance. We expect very similar methods to give similar results if loops are permitted.

\subsection{Contributions}

Our main contributions are giving a complete characterization of eventual forever connectedness in these models and, in this connected regime, proving Janson's \cite{janson_2017_on} conjecture on the asymptotic normality of the number of vertices. We also resolve an open problem of Janson \cite{janson_2017_on} regarding when the produced graphs are eventually forever essentially complete and obtain an apparently new result characterizing when the urns in a classical urn scheme \cite{benhamou_2017_concentration,hwang_2008_local,dutko_1989_central,karlin_1967_central,gnedin_2007_notes} are eventually forever filled in order.

\subsection{Structure} The remainder of Section \ref{Introduction} consists of two parts. The first, Subsection \ref{Motivation}, will motivate the model of interest and our contribution to the theory thereof. The second, Subsection \ref{Related work}, will discuss previous work on the model we are studying, as well as on other related models. In Section \ref{Preliminaries and Notation} we introduce some preliminaries and notation. Each of the following three sections, namely Sections \ref{Connectedness}, \ref{Gaussianity} and \ref{Completeness}, is dedicated to one of the three main theorems. The first section covers the characterization of eventual forever connectedness, the second the asymptotic Gaussianity of the vertex count in the connected regime and the third the characterization of eventual forever essential completeness.

\subsection{Motivation} \label{Motivation}

\subsubsection{Motivating Edge Exchangeable Models}
A sequence of random variables is said to be exchangeable if its law is invariant to finite permutations. This is weaker than being i.i.d. from some distribution, instead capturing the intuition that the random variables are observed in no particular order. The most classical notion of exchangeability for random graphs is vertex exchangeability. For example, if the network is the friendships on a social media site, then vertex exchangeability would mean assuming that users join the website in no particular order. However, Crane and Dempsey, paraphrasing and summarizing previous authors, on the topic of vertex exchangeability, write "A realization from an exchangeable model for countable graphs is
sparse if and only if it has no edges" \cite{crane_2018_edge}. Since graphs observed in practice are typically sparse, vertex exchangeable graphs are therefore not satisfactory. In many applications it is natural to instead assume that interactions (edges) are observed in no particular order \cite{crane_2018_edge,cai_2016_edgeexchangeable}, that is, to assume edge exchangeability. The models we study are edge exchangeable but are not quite the most general edge exchangeable graphs possible, principally because $\mu$ could in general itself be random. See \cite{crane_2018_edge} for a general characterization theorem.

\subsubsection{Motivating Our Contribution}  For applications in statistics one would like to ensure that edge exchangeable models, unlike vertex exchangeable models, can produce graphs which resemble those encountered in practice. This calls for understanding how the qualities of the produced graphs depend on $\mu$. Moreover, in some applications one may have a choice of (prior over) $\mu$. In that case results of the type we provide may help inform the choice of this infinite-dimensional parameter. Moving beyond motivations in bona fide statistical network theory, the models, and in particular their hypergraph generalizations, appear in the generalized unseen species problem \cite{eriksson_2026_the}. Theorem \ref{connected Gauss} would then, in the phone calls application discussed in \cite{eriksson_2026_the}, establish that if one is spying on a terrorist organization by monitoring their phone calls, then the number of discovered terrorists in a finite time interval has, under suitable assumptions, approximately Gaussian fluctuations. This is a stepping stone towards constructing calibrated confidence intervals for linear estimators in the generalized unseen species problem.

In addition to being motivated by statistical applications, the model is quite natural and displays some interesting phenomena which are not typical of random graph models, making them probabilistically interesting. For example, the vertex set is not deterministic, we prove results which are really about the entire history of a graph-valued process, not just its marginals and we find that a topological assumption has distributional consequences.

\subsection{Related Work} \label{Related work}

\subsubsection{Previous Work on Edge Exchangeable Models} Edge exchangeable models go back to at least \cite{crane_2018_edge,cai_2016_edgeexchangeable}. Both \cite{crane_2018_edge} and \cite{cai_2016_edgeexchangeable} consider sparsity both before and after identifying parallel edges. They both show that sparsity is possible for edge exchangeable graphs by showing that a particular model is sparse. A more general study of sparsity as well as other graph properties was undertaken by Janson in \cite{janson_2017_on}. The present paper can be seen as a continuation of the research program started in \cite{janson_2017_on}. Janson gave formulas for the expected number of distinct edges and vertices and showed these quantities to be strongly relatively stable. Moreover, asymptotic normality for the number of (distinct) edges was established under the assumption that the variance of the number of edges goes to infinity. Example $\mu$ are given which produce everything from extremely sparse to extremely dense graphs. A sufficient condition on $\mu$ is given such that the graphs are eventually forever essentially complete (see Definition \ref{essentially complete definition}). Statistical applications of edge exchangeability are also being developed, for instance, anomaly detection \cite{luo_2023_anomalous}. The theory of simulation (sampling) and inference for a class of $\mu$ was developed in \cite{li_2021_truncated}. The textbook by Crane \cite{crane_2018_probabilistic} contains an introduction to edge exchangeable models.

\subsubsection{Comparisons to other Models} \label{Comparisons to other Models}
Another type of exchangeability for graph-valued stochastic processes is that of Caron and Fox's exchangeable measures \cite{caron_2017_sparse}. These models can also produce sparse graphs, but as Crane expresses in his comment in \cite{caron_2017_sparse}, it is much less clear what their exchangeability corresponds to in terms of modeling the real world and hence when it can be expected.

Models which resemble the ones we are interested in arise in various contexts. For instance, we may imagine an inference problem in which there is some true unknown finite graph and for each edge there is an exponential time (independent but not necessarily identically distributed) until we observe it. This was considered, and suitable statistical estimators developed, in \cite{aldous_2021_a}. The same model, although therein viewed from the perspective of percolation theory more so than statistics, is considered in \cite{aldous_2016_weak,aldous_2019_the} where some concentration results for the emergence of a giant component are established. This has some rough parallels to our characterization of connectedness.

One general type of graph model is the so-called "inhomogeneous random graph". There one works with a finite deterministic vertex set and includes edge $e$, independently of all other edges, with some probability $p_e$. At all times our exponential arrival time graph can be viewed as an instantiation of this model with an infinite vertex set, $p_e = 1-e^{-\mu_e t}$ and the additional rule that all isolated vertices are discarded.  

In this work we are agnostic as to how $\mu$ and consequently the $p_e$ arise. However, every way to obtain $p_e$ for an inhomogeneous random graph can also produce $\mu_e$. For instance, one may imagine that each vertex is associated to a vector in some latent space and that the $\mu_e$ arise as dot products of these vectors like the classical dot product random graphs, a survey of which can be found in \cite{athreya_2017_statistical}. Note that in this case some technical care needs to be taken to ensure that $\mu$ can be normalized.

The Norros–Reittu model \cite{norros_2006_on} is another model which has an independent Poisson number of edges between each pair of vertices. However, in that model the intensity is (up to a common scaling) the product of intensities associated to the vertices. Moreover, the vertex set is deterministic unlike the vertex set in our model. 

\section*{Acknowledgments}
We thank Alexander Fuchs-Kreiss for many useful conversations regarding these ideas and for offering valuable feedback on earlier drafts.

\section{Preliminaries and Notation} \label{Preliminaries and Notation}
For a set $\mathcal{X}$ we use $\mathcal{P}(\mathcal{X})$ to denote the set of all probability measures over $\mathcal{X}$. When clear from context we drop limits of summation. For $\mu \in \mathcal{P}(\mathbb{N}_2)$ define $\mu_{ij}:=\mu(\{i,j\})$ and define $M_i = \sum_{j}\mu_{ij}$ for $i \in \mathbb{N}$. For $e=\{i,j\}$ define $M_e := M_{ij} := M_i+M_j-\mu_{ij}$.

\begin{definition}
    For a stochastic process $X$ indexed by $I$ taking values in some set $\mathcal{X}$ and a property $A: \mathcal{X} \rightarrow \{\text{True},\text{False}\}$ we say that $X$ has or satisfies $A$ \textbf{eventually forever} if and only if there exists an $i$ such that $A(X(j)) = \text{True}$ for all $j \in I$ such that $j>i$. Equivalently, if the indicators of that property converge to $1$.
\end{definition}
For the processes and properties we are interested in, whether the process satisfies the property eventually forever has probability either zero or one. We may therefore omit the "almost surely" and simply write $X$ has $A$ eventually forever to mean that $X$ almost surely has $A$ eventually forever.
\begin{lemma}\label{De-Poissonization Lemma} If $A$ is a property that depends only on the order of appearance of edges, then $G_t$ has $A$ eventually forever if and only if $G_n$ has $A$ eventually forever.
\end{lemma} \begin{proof}
    This is immediate from the fact that we can sample $G_n$ by waiting until $n$ events have occurred in $G_t$.
\end{proof}

For a family of random variables indexed by $I$ and for nonempty $J \subset I$ we write $\check{J}$ and $\hat{J}$ for the supremum and infimum among random variables indexed by elements of $J$. For our applications these are always attained, so one may speak of the maximum and the minimum, respectively.

\begin{definition}\label{asymptotic gauss definition} Let $X_t$ be a stochastic process in either discrete or continuous time which takes values in (a subset of) $\mathbb{R}$. We call such a process \textbf{asymptotically Gaussian} if there exist deterministic functions $f,g: \mathbb{R} \to \mathbb{R}$ such that
    \begin{align*}
        \frac{X_t-f(t)}{g(t)} \rightarrow \mathcal{N}(0,1),
    \end{align*} where $\mathcal{N}(0,1)$ is the standard normal distribution and the convergence is in distribution as $t \rightarrow \infty$.
\end{definition}
\section{Connectedness}
\label{Connectedness}

We will now build up towards characterizing eventual forever connectedness. 

\begin{definition}
    For a $\mu \in \mathcal{P}(\mathbb{N}_2)$ we call the graph $(\cup_{e: \mu_e>0}e,\{e \in \mathbb{N}_2: \mu_e>0\})$ the \textbf{support} of $\mu$. When it is a connected graph we say that $\mu$ has \textbf{connected support}.
\end{definition}
Let $I_{e}$ be the event that the edge $e$ brings two new vertices to $G_t$ when it arrives. This is well-defined since almost surely edges are not added at the same time. 

\begin{lemma} \label{fin components imply eventually connected} Outside of an event of probability zero, $G_t$ is eventually forever connected if and only if only finitely many $I_e$ occur and $\mu$ has connected support.
\end{lemma}
\begin{proof}
First we show that the violation of either condition means that eventual forever connectedness fails. If $\mu$ does not have connected support, then $G_t$ is not eventually forever connected since edges that can never belong to the same component will eventually appear. When an edge arrives it adds zero, one or two new vertices. A new component is created if and only if the edge brings two new vertices. Thus if infinitely many $I_e$ occur there are infinitely many times at which the graph is not connected. Since a.s. only finitely many edges are added in finite time this means that the graph will not be eventually forever connected.

Next we show that if both conditions are satisfied then eventual forever connectedness follows. Since only finitely many $I_e$ occur, eventually all of them will have occurred. Since $\mu$ has connected support, edges which cause all the created components to join up will subsequently appear, after which point $G_t$ will be connected.
\end{proof} 

Lemma \ref{fin components imply eventually connected} reduces the problem of determining eventual forever connectedness to deciding whether infinitely many $I_e$ occur. We do this by computing $P(I_e)$ and $P(I_e \wedge I_f)$ and applying the first Borel-Cantelli Lemma and a suitable partial converse. We now define the notation that will be required for the computation of $P(I_e \wedge I_f)$. For $e,f \in \mathbb{N}_2$ with $e \cap f = \emptyset$ let $C_{ef} = \{g \in \mathbb{N}_2: g \cap e \neq \emptyset, g \cap f \neq \emptyset\}$, $C_e := \{g: |e \cap g|=1\}\backslash C_{ef}$ and define $C_f$ analogously. Note that $\{e\},\{f\},C_{ef},C_e,C_f$ are disjoint and that $\tau_e,\tau_f,\hat{C}_{ef},\hat{C}_e,\hat{C}_f$ are therefore independent and exponentially distributed with parameters $\mu_e,\mu_f,b_{ef} := \sum_{g \in C_{ef}}\mu_g,r_e := \sum_{g \in C_{e}}\mu_g,r_f := \sum_{g \in C_{f}}\mu_g$ respectively. For compactness of notation we define $a_e := r_e+\mu_e$, $c_e := 1-\frac{r_e}{r_e+\mu_e}=\frac{\mu_e}{a_e}$. We include Figure \ref{fig:notation} to help with parsing this notation.

\begin{figure}[t]
\centering
\begin{tikzpicture}
\node[fill=black, circle, inner sep=1.5pt] (e1) at (-2,2) {};
\node[fill=black, circle, inner sep=1.5pt] (e2) at (2,2) {};
\node[fill=black, circle, inner sep=1.5pt] (f1) at (-2,0) {};
\node[fill=black, circle, inner sep=1.5pt] (f2) at (2,0) {};
\node[fill=black, circle, inner sep=1.5pt] (ex1) at (-4,2.5) {};
\node[fill=black, circle, inner sep=1.5pt] (ey1) at (4,2.5) {};
\node[fill=black, circle, inner sep=1.5pt] (fx1) at (-4,1) {};
\node[fill=black, circle, inner sep=1.5pt] (fx2) at (-4,-0.5) {};
\node[fill=black, circle, inner sep=1.5pt] (fy1) at (4,1) {};
\node[fill=black, circle, inner sep=1.5pt] (fy2) at (4,-0.5) {};

\draw[thick] (e1) -- (e2) node[midway,above] {$e$};
\draw[thick] (f1) -- (f2) node[midway,above] {$f$};
\draw[thick] (f1) -- (fx1) node[midway,above] {$C_f$};
\draw[thick] (f1) -- (fx2) node[midway,above] {$C_f$};
\draw[thick] (f2) -- (fy1) node[midway,above] {$C_f$};
\draw[thick] (f2) -- (fy2) node[midway,above] {$C_f$};
\draw[thick] (e1) -- (ex1) node[midway,above] {$C_{e}$};
\draw[thick] (e1) -- (fx1) node[midway,above] {$C_{e}$};
\draw[thick] (e2) -- (ey1) node[midway,above] {$C_{e}$};
\draw[thick] (e2) -- (fy1) node[midway,above] {$C_{e}$};

\draw[thick] (e1) -- (f1) node[midway,right] {$C_{ef}$};
\draw[thick] (e1) -- (f2) node[midway,above] {$C_{ef}$};

\end{tikzpicture}
\caption{Some selected edges labeled by the sets they belong to}
\label{fig:notation}

\end{figure}

\begin{lemma}\label{expl prob} For $e,f \in \mathbb{N}_2$ we have $|e \cap f| \in \{0,1,2\}$ and
\begin{itemize} 
    \item $|e \cap f| = 0 \implies P(I_e \wedge I_f) = c_ec_f\left(1-\frac{b_{ef}}{a_e+b_{ef}}-\frac{b_{ef}}{a_f+b_{ef}}+\frac{b_{ef}}{a_e+a_f+b_{ef}}\right)$.
    \item $|e \cap f| = 1 \implies P(I_e \wedge I_f) = 0$.
    \item $|e \cap f| = 2 \implies e = f \text{ and } P(I_e \wedge I_f)=\frac{\mu_e}{M_e}$ if $\mu_e > 0$ and $P(I_e)=0$ otherwise.
\end{itemize}
\end{lemma} \begin{proof} $|e \cap f| \in \{0,1,2\}$ is immediate from the fact that $|e|=|f|=2$. If either $\mu_e$ or $\mu_f$ is zero, it is clear that $P(I_e \wedge I_f) = 0$ in which case the lemma holds, so we may assume both to be non-zero. First we deal with the case $|e \cap f| = 0$. The case $b_{ef} =0$ is trivial so we may assume $b_{ef} > 0$. We will need the following sub-lemma.
\begin{align*}
P(I_e|\hat{C}_{ef}=t) &= \int_0^\infty P(I_e|\hat{C}_{ef}=t,\hat{C}_e = \tau)r_ee^{-r_e \tau}d\tau, \\
&= \int_0^\infty (1-e^{-\mu_e\min(t,\tau)}) r_ee^{-r_e \tau}d\tau, \\
&= 1-\int_0^\infty e^{-\mu_e\min(t,\tau)} r_ee^{-r_e \tau}d\tau, \\
&= 1-\int_0^te^{-\mu_e \tau} r_e e^{-r_e \tau}d\tau-e^{-\mu_e t}\int_t^\infty r_ee^{-r_e \tau}d\tau,\\
&= 1-\frac{r_e}{r_e+\mu_e}\left(1-e^{-(r_e+\mu_e)t}\right)-e^{-(\mu_e+r_e)t},\\
&= \left(1-\frac{r_e}{r_e+\mu_e}\right)\left(1-e^{-(r_e+\mu_e)t}\right), \\
&= c_e(1-e^{-a_e t}).
\end{align*}
Now we compute
    \begin{align*}
        P(I_e \wedge I_f) &= \int_0^\infty P(I_e \wedge I_f|\hat{C}_{ef}=t)b_{ef}e^{-b_{ef}t}dt, \\
        &=  \int_0^\infty P(I_e|\hat{C}_{ef}=t)P(I_f|\hat{C}_{ef}=t)b_{ef}e^{-b_{ef}t}dt, \\
        &=  \int_0^\infty c_e(1-e^{-a_e t})c_f(1-e^{-a_f t})b_{ef}e^{-b_{ef}t}dt, \\
        &= c_ec_f\left(1-\frac{b_{ef}}{a_e+b_{ef}}-\frac{b_{ef}}{a_f+b_{ef}}+\frac{b_{ef}}{a_e+a_f+b_{ef}}\right).
    \end{align*} where the third equality is obtained by applying the sub-lemma twice.
    
If $|e \cap f|=1$ let $ v \in e \cap f$ be the shared vertex. For both $I_e$ and $I_f$ to occur both edges have to bring $v$. Since $e \neq f$, edges almost surely do not arrive at the same time and the vertex set is increasing, this is not possible, so $P(I_e \wedge I_f) = 0$.

If $|e \cap f| = 2$, then, since $|e|=|f|=2$ we have $e=f$ and $I_e \wedge I_f = I_e$.
The event $I_e$ is precisely the event that $\hat{\{g: e \cap g \neq \emptyset \}} = \tau_e$ and so, by a standard property of exponential random variables,
    $$P(I_e)=\frac{\mu_e}{M_e}.$$ 
\end{proof}
\begin{theorem} \label{Eventual forever connectedness thm}
    $G_t$ and $G_n$ are eventually forever connected with probability one if $\sum_{e \in \mathbb{N}_2}\frac{\mu_e}{M_e} < \infty$ and $\mu$ has connected support. Else they are almost surely not eventually forever connected.
\end{theorem} \begin{proof} By Lemma \ref{De-Poissonization Lemma} it suffices to consider $G_t$. By Lemma \ref{fin components imply eventually connected}  it suffices to resolve whether infinitely many $I_e$ occur. If $$\sum_{e \in \mathbb{N}_2} \frac{\mu_e}{M_e} < \infty,$$ then almost surely only finitely many $I_e$ occur by the first Borel-Cantelli lemma and the third case of Lemma \ref{expl prob}. Now assume instead that $\sum_{e \in \mathbb{N}_2} \frac{\mu_e}{M_e} = \infty$.  For $|e \cap f| = 0$ we have,

    \begin{align*}
        \frac{P(I_e \wedge I_f)}{P(I_e)P(I_f)} &= \frac{a_e+b_{ef}}{\mu_e}\frac{a_f+b_{ef}}{\mu_f}c_ec_f\left(1-\frac{b_{ef}}{a_e+b_{ef}}-\frac{b_{ef}}{a_f+b_{ef}}+\frac{b_{ef}}{a_e+a_f+b_{ef}}\right), \\
        &= \frac{a_e+a_f+2b_{ef}}{a_e+a_f+b_{ef}}, \\
        &< 2.
    \end{align*} We used the first and third cases in Lemma \ref{expl prob} to obtain the first equality, along with the observation that $M_e = a_e+b_{ef}$. The second equality was obtained by lengthy but elementary algebra. The inequality is immediate from noticing that the numerator is less than twice the denominator.
    
    Since $\mathbb{N}_2$ is countable we may order it in some arbitrary but fixed way. In the following we take $e < n$ to mean that $e \in \mathbb{N}_2$ is one of the first $n-1$ edges under the ordering.
\begin{align*}
    \sum_{e,f < n}P(I_e \wedge I_f) = \sum_{e < n}P(I_e) + \sum_{e,f < n: e \neq f}P(I_e \wedge I_f).
\end{align*}
Denoting the first sum by $A_n$ and the second by $B_n$,

\begin{align*}
    B_n \leq 2\sum_{e,f < n: e \neq f} P(I_e)P(I_f) \leq 2A_n^2. 
\end{align*}
Therefore

\begin{align*}
    \frac{A_n^2}{A_n + B_n} \geq \frac{A_n^2}{A_n+2A_n^2} \rightarrow \frac{1}{2},
\end{align*}
since $A_n \rightarrow \infty$ by assumption.
Applying the Kochen-Stone Lemma \cite{yan_2006_a} now gives $P(I_e \text{ i.o.}) \geq \frac{1}{2}$.

We will next argue that $I_e$ \text{ i.o.} is a tail event in the sense of Kolmogorov's 0-1 law. Because each vertex can be involved in only one double vertex event, changing an edge arrival time can change $\sum_e 1_{I_e}$ by at most $2$. More explicitly, consider making the edge $e$ arrive sooner. This could create one more double-new-vertex-event, namely $I_e$. It could also destroy up to two such events, at most one for each vertex in $e$. On net it changes the number of such events by at most $2$. By symmetry this also covers increasing $\tau_e$. From this one sees that changing the arrival times in a finite collection can only change the number of double vertex events by a finite amount, namely by twice the size of the collection. This cannot change whether infinitely many $I_e$ occur. Moreover, the edge arrivals are independent, so indeed whether infinitely many $I_e$ occur is a tail event. Invoking Kolmogorov's 0-1 law now completes the proof.
\end{proof}
\begin{example} \label{gamma > 2} Let $\mu_{ij} \propto (ij)^{-\gamma}$ for $i \neq j$ and some $\gamma > 1$. Then $G_n$ and $G_t$ are eventually forever connected if and only if $\gamma > 2$.
\end{example}\begin{proof} Note that we need $\gamma > 1$ for $\mu$ to be normalizable. The support of $\mu$ is the complete graph which is connected. Note that $\frac{\mu_{e}}{M_{e}}$ is invariant under rescaling of $\mu$ so we may assume that the proportionality is in fact an equality. We compute
\begin{align*}
    M_{ij} &= \sum_{k}\mu_{kj}+\sum_k\mu_{ik}-\mu_{ij} \\
    &=\sum_{k \neq j}(kj)^{-\gamma}+\sum_{k \neq i}(ik)^{-\gamma}-(ij)^{-\gamma} \\
    &= j^{-\gamma}\zeta(\gamma)-(jj)^{-\gamma}+i^{-\gamma}\zeta(\gamma)-(ii)^{-\gamma}-(ij)^{-\gamma},
\end{align*} where $\zeta$ denotes the Riemann zeta function. This implies that 
\begin{align} \label{gamma double sum}
    \sum_{ij}\frac{\mu_{ij}}{M_{ij}} &= \sum_{i \neq j}\frac{(ij)^{-\gamma}}{j^{-\gamma}\zeta(\gamma)-(jj)^{-\gamma}+i^{-\gamma}\zeta(\gamma)-(ii)^{-\gamma}-(ij)^{-\gamma}}, \\
    &= \sum_{i \neq j}\frac{1}{i^\gamma(\zeta(\gamma)-j^{-\gamma})+j^\gamma(\zeta(\gamma)-i^{-\gamma})-1}.
\end{align} 
Notice that 
$$\lim_{i,j \rightarrow \infty} \frac{i^\gamma(\zeta(\gamma)-j^{-\gamma})+j^\gamma(\zeta(\gamma)-i^{-\gamma})-1}{i^\gamma+j^\gamma} = \zeta(\gamma) > 0.$$ 
From the limit comparison test for double sums, Theorem 9.6 of \cite{mursaleen_2013_double}, it follows, after some technicalities, that $\sum_e \frac{\mu_e}{M_e}$ converges if and only if $$\sum_{i=1}^\infty \sum_{j=1}^\infty \frac{1}{i^\gamma+j^\gamma},$$ converges. Finally, by Propositions 1 and 2 of \cite{zanu_2025_on} this sum converges if and only if $\gamma > 2$.
\end{proof}
\begin{remark}
    From \cite{janson_2017_on} it is known that $\mu$ of this type produce sparse graphs regardless of $\gamma$. We thus see that edge exchangeable graphs can simultaneously have good connectedness and good sparsity properties.
\end{remark}
Since a.s. convergence is stronger than convergence in probability we have the following.
\begin{corollary}\label{marginal pass} If $\mu$ has connected support and $$\sum_e \frac{\mu_e}{M_e} < \infty,$$ then the probability that the graphs are connected goes to one.
\end{corollary}

\section{Gaussianity} \label{Gaussianity}
Janson \cite{janson_2017_on} showed asymptotic Gaussianity for the number of edges by coupling to certain urn schemes and we take a similar approach for the vertices. The key difference is that the arrivals of edges are, after Poissonizing, independent while the arrivals of the vertices are not. We will assume that the graph is eventually forever connected. Our key insight is that this assumption will allow us to control the dependence.

We will now introduce two urn schemes. For the first scheme, fix a sequence of probabilities $\{p_i\}_{i=1}^\infty$ such that $\sum_i p_i = 1$. For each $i \in \mathbb{N}$ define an urn. At times $n=1,2,3,\ldots$ throw a ball with probability $p_i$ of landing in urn $i$. Let $U_n$ be the set of occupied urns at time $n$ so that $|U_n|$ is the number of occupied urns at time $n$. For the second scheme, fix a sequence of non-negative reals $\{\lambda_i\}_{i=1}^\infty$ such that $\sum_i \lambda_i < \infty$. For each urn define a Poisson process with intensity $\lambda_i$ and add a ball to the urn whenever there is an event in the associated Poisson process. Let $U_t$ be the set of occupied urns at time $t$ so that $|U_t|$ is the number of occupied urns at time $t$. When $\sum_i \lambda_i = 1$ the latter scheme is the Poissonized version of the former with $p_i = \lambda_i$. We will abuse notation and write $U_t$ even when we are evaluating the continuous time scheme at the integer time $n$, in this sense $t$ and $n$ should be thought of as interchangeable and which one we use merely distinguishes continuous and discrete time. Convergence to $\mathcal{N}(0,1)$ will always be in distribution as time goes to infinity.

\begin{theorem} (\cite{dutko_1989_central}, Theorems 1 and 2) \label{Urn Gauss}
    Let $U_n, U_t$ be as above. If $\mathbb{V}[|U_t|] \rightarrow \infty$, then
    \begin{align*}
        \frac{|U_t|-\mathbb{E}[|U_t|]}{\sqrt{\mathbb{V}[|U_t|]}} \rightarrow \mathcal{N}(0,1), \\
        \frac{|U_n|-\mathbb{E}[|U_n|]}{\sqrt{\mathbb{V}[|U_t|]}} \rightarrow \mathcal{N}(0,1).
    \end{align*}
\end{theorem}
To be precise, \cite{dutko_1989_central} assumed $\sum_j \lambda_j =1$ but the case $\sum_j \lambda_j \neq 1$ is handled by a simple scaling argument. Asymptotic Gaussianity of the number of (distinct) edges is a corollary of Theorem \ref{Urn Gauss}, obtained by defining an urn for each edge \cite{janson_2017_on}. For the number of vertices it is an open problem, namely Problem 6.9 of \cite{janson_2017_on}. While Janson only conjectures the result in discrete time we will also consider continuous time.

The indicators of the vertices, unlike the edges, are dependent. However, we will show that actually one can couple to the standard urn scheme, if one assumes eventual forever connectedness.

\begin{theorem} \label{connected Gauss}
    Let $G_t$ be eventually forever connected and suppose that $\mathbb{V}[|U_t|] \rightarrow \infty$. Then the following hold. \begin{align*}
    \frac{|V_n|-\mathbb{E}[|V_n|]}{\sqrt{\mathbb{V}[|U_t|]}} \rightarrow \mathcal{N}(0,1), \\
        \frac{|V_t|-\mathbb{E}[|V_t|]}{\sqrt{\mathbb{V}[|U_t|]}} \rightarrow \mathcal{N}(0,1).
    \end{align*}
\end{theorem}
We begin by sketching the proof. We will couple to the urn scheme with intensities $\lambda_j = M_j$. By choosing an appropriate coupling, vertices will be added to $V_t$ but not $U_t$ only when a pair of new vertices are simultaneously added to $V_t$. Since we know this happens only finitely many times in the connected regime, this will not cause us any trouble. We will need to fill some extra urns in $U_t$ but no more often than the double vertex additions to $V_t$. Thus there will only be finitely many vertices that are ever in one but not the other. Eventually these vertices will be present in both and from that point onward $U_t$ and $V_t$ will be equal. The coupling inequality will then complete the proof.
\begin{proof} We define a coupling $(\tilde{V}_t,\tilde{U}_t)$ of $V_t$ and $U_t$. Start with $\tilde{V}_t=\tilde{U}_t =\emptyset$.
Repeatedly and independently wait an $\exp(3)$ amount of time, then do the following.

For every $i \in \mathbb{N} \backslash \tilde{U}_t$ compute 
\[  \lambda^{(i)}(\tilde{V}_t,\tilde{U}_t) := \left\{
\begin{array}{ll}
      \frac{1}{2}\sum_{j \in \overline{\tilde{V}}_t \cap \overline{\tilde{U}}_t}\mu_{ij} & \text{if } i \notin \tilde{V}_t,\\
      \frac{1}{2}\sum_{j \in \overline{\tilde{U}}_t \cap \tilde{V}_t}\mu_{ij} + \sum_{j \in \overline{\tilde{U}}_t \cap \overline{\tilde{V}}_t}\mu_{ij} & \text{if } i \in \tilde{V}_t. \\
\end{array} 
\right. \]
Then generate a $\mathbb{N} \cup \mathbb{N}_2 \cup \{\emptyset\}$ -valued random variable $\chi$ by generating $\eta \sim \text{Uniform}[0,1]$ independent of everything else and thresholding appropriately to obtain:
\begin{itemize}
    \item $P(\chi = i) = \frac{\lambda^{(i)}}{3}$ for all $ i \in \mathbb{N}$.
    \item $P(\chi = e) = \frac{\mu_{e}}{3}$ for all $e \in \mathbb{N}_2$.
    \item $P(\chi = \emptyset) = \frac{2}{3}-\frac{1}{3}\sum_i \lambda^{(i)}.$ 
\end{itemize}
Then increment the processes according to
\begin{itemize}
    \item If $\chi = \emptyset$: pass 
    \item \textcolor{blue}{If $\chi = i$: Add $i$ to $\tilde{U}_t$}.
    \item If $\chi = \{i,j\}$ for some $\{i,j\} \in \mathbb{N}_2$, 
    \begin{itemize}
        \item If $|\{i,j\}\backslash \tilde{V}_t|=0$:
        \begin{itemize}
            \item If $|\{i,j\}\backslash \tilde{U}_t| = 0$: pass.
            \item \textcolor{green}{If $|\{i,j\} \backslash \tilde{U}_t| =1$: add $\{i,j\} \backslash \tilde{U}_t$ to $\tilde{U}_t$}.
            \item \textcolor{magenta}{If $|\{i,j\} \backslash \tilde{U}_t| = 2$: Choose $i$ or $j$ according to a coin toss, independent of everything else, and add that one to $\tilde{U}_t$.}
        \end{itemize}
        \item If $|\{i,j\}\backslash \tilde{V}_t|=1$:
        Let $x \in \{i,j\}\backslash \tilde{V}_t$ be that one element and let $y \in \{i,j\} \cap \tilde{V}_t$ be the other. Then:
        \begin{itemize}
        \item \textcolor{yellow}{If $x \notin \tilde{U}_t$: Add $x$ to $\tilde{U}_t$}.
        \item \textcolor{violet}{If $x \in \tilde{U}_t$ and $y \notin \tilde{U}_t$: add $y$ to $\tilde{U}_t$}.
        \item If $x \in \tilde{U}_t$ and $y \in \tilde{U}_t$: pass.
        \end{itemize}
        
        \item If $|\{i,j\}\backslash \tilde{V}_t|=2$:
        \begin{itemize}
            \item If $|\{i,j\}\backslash \tilde{U}_t| = 0$: pass.
            \item \textcolor{orange}{If $|\{i,j\} \backslash \tilde{U}_t| =1$: add $\{i,j\} \backslash \tilde{U}_t$ to $\tilde{U}_t$}.
            \item \textcolor{brown}{If $|\{i,j\} \backslash \tilde{U}_t| = 2$: Choose $i$ or $j$ according to a coin toss, independent of everything else, and add that one to $\tilde{U}_t$.}
        \end{itemize}
        
        \item Add $\{i,j\} \backslash \tilde{V}_t$ to $\tilde{V}_t$
    \end{itemize}
\end{itemize}
To ensure that the above construction is well-defined we need to ensure that $\chi$ is well-defined. To see this note that $\lambda^{(i)} \leq \sum_{j}\mu_{ij} = M_i$ and that $ \sum_i M_i = 2,$ so that therefore $$P(\chi = \emptyset) = \frac{2}{3}-\frac{1}{3}\sum_i\lambda^{(i)} \geq 0.$$ Note also that $$\sum_iP(\chi = i)+\sum_e P(\chi = e) + P(\chi=\emptyset) = \sum_i\frac{\lambda^{(i)}}{3}+\frac{1}{3}+\frac{2}{3}-\sum_i\frac{\lambda^{(i)}}{3}= 1.$$

From the thinning property of Poisson processes and the third bullet point in the increment rule we see that $\tilde{V}_t$ has the same distribution as $V_t$.

Let $p$ be the probability that $i \notin \tilde{U}_t$ is added to $\tilde{U}_t$ after a particular $\exp(3)$ wait. By summing over the probabilities of all the different ways $i$ can be added to $\tilde{U}_t$ this can be seen to be given by

\begin{align*}
    3p = &\textcolor{blue}{1_{i \notin \tilde{V}_t}\frac{1}{2}\sum_{j \in \overline{\tilde{V}}_t \cap \overline{\tilde{U}}_t}\mu_{ij}+1_{i \in \tilde{V}_t}\frac{1}{2}\sum_{j \in \tilde{V}_t \cap \overline{\tilde{U}}_t}\mu_{ij}+1_{i \in \tilde{V}_t}\sum_{j \in \overline{\tilde{V}}_t \cap \overline{\tilde{U}}_t}\mu_{ij}}, \\
    +&\textcolor{magenta}{1_{i \in \tilde{V}_t}\frac{1}{2}\sum_{j \in \tilde{V}_t \cap \overline{\tilde{U}}_t}\mu_{ij}} + \textcolor{green}{1_{i \in \tilde{V}_t}\sum_{j \in \tilde{V}_t \cap \tilde{U}_t}\mu_{ij}}+\textcolor{yellow}{1_{i\notin \tilde{V}_t}\sum_{j \in \tilde{V}_t}\mu_{ij}}\\
    +&\textcolor{violet}{1_{i \in \tilde{V}_t} \sum_{j \in \overline{\tilde{V}}_t \cap \tilde{U}_t}\mu_{ij}} + \textcolor{orange}{1_{i \notin \tilde{V}_t}\sum_{j \in \overline{\tilde{V}}_t \cap \tilde{U}_t}\mu_{ij}}
    +\textcolor{brown}{1_{i \notin \tilde{V}_t}\frac{1}{2}\sum_{j \in \overline{\tilde{V}}_t \cap \overline{\tilde{U}}_t}\mu_{ij}} \\
    = &1_{i \in \tilde{V}_t}\sum_{j}\mu_{ij}+1_{i \notin \tilde{V}_t}\sum_{j}\mu_{ij}, \\
    = &M_i.
\end{align*} Note that this does not depend on $t,\tilde{V}_t$ and $\tilde{U}_t$. Moreover, there are no simultaneous jumps, the $\eta$ are independent, as are the $\exp(3)$ waiting times, so looking only at $\tilde{U}_t$ we observe an urn filling process and therefore $(|\tilde{V}_t|,|\tilde{U}_t|)$ is a valid coupling of $|V_t|,|U_t|$.

Next we show that $\tilde{V}_t = \tilde{U}_t$ eventually forever by showing eventual forever containment both ways. We begin by showing that $\tilde{V}_t \subset \tilde{U}_t$ eventually forever. First, notice that due to how we have set up the coupling the only way for a vertex $v$ to be added to $\tilde{V}_t \backslash \tilde{U}_t$ is by an edge $\chi$ adding both $v$ and another vertex to $\tilde{V}_t$. Since this happens only finitely many times, by the connectedness assumption, eventually, it will cease to happen and, since those vertices will eventually appear in $\tilde{U}_t$, they will no longer be in $\tilde{V}_t \backslash \tilde{U}_t$. Thus, eventually forever $\tilde{V}_t \subset  \tilde{U}_t$. Next we show $\tilde{U}_t \subset \tilde{V}_t$ eventually forever. Immediately following an $\exp(3)$ wait, the probability that an urn is filled in $\tilde{U}_t$ without the corresponding vertex being added to or already being in  $\tilde{V}_t$ is $\frac{1}{6}\sum_{u \notin \tilde{V}_t} \sum_{v \in \overline{\tilde{U}}_t \cap \overline{\tilde{V}}_t}\mu_{uv}$. This is because we need a $\textcolor{blue}{\chi = i}$ event where $i \notin \tilde{V}_t$. Now observe that after one of the $\exp(3)$ waits the probability of adding two new vertices to $\tilde{V}_t$ at the same time is $\frac{1}{3}\sum_{i,j \in \overline{\tilde{V}}_t}\mu_{ij}$ and that

\begin{align}\label{intrinsic dom}\frac{1}{6}\sum_{i \notin \tilde{V}_t} \sum_{j \in \overline{\tilde{U}}_t \cap \overline{\tilde{V}}_t}\mu_{ij} \leq \frac{1}{3}\sum_{i,j \in \overline{\tilde{V}}_t}\mu_{ij}.\end{align}

Since the $\eta$ are independent, (\ref{intrinsic dom}) is sufficient to see that the distribution of the number of positive integers that are ever in $\tilde{U}_t \backslash \tilde{V}_t$ is dominated by the number of vertices that appear at the same time as another vertex in $\tilde{V}_t$. Since the latter is a.s. finite, so is the former. This is perhaps easier to see in terms of processes, in which case it corresponds to the obvious fact that if a Cox-process only ever has finitely many events and another process has a lower intensity, then so does the second. Since any vertices added to $\tilde{U}_t \backslash \tilde{V}_t$ eventually leave, one has $\tilde{U}_t \backslash \tilde{V}_t = \emptyset$ eventually forever. Thus we have two-way containment and $\tilde{V}_t = \tilde{U}_t$ eventually forever. In particular $|\tilde{V}_t|=|\tilde{U}_t|$ eventually forever almost surely.

By the coupling inequality and since almost sure convergence is stronger than convergence in probability it follows that 
\begin{align}\label{TV vertex urn} d_{TV}(\mathcal{L}[|V_t|],\mathcal{L}[|U_t|]) \leq P(|\tilde{V}_t| \neq |\tilde{U}_t|) \rightarrow 0,\end{align} where $\mathcal{L}$ denotes the law of a random variable. Now combining the fact that the Lévy–Prokhorov distance, which we denote by $\pi$, metrizes convergence in distribution \cite{gibbs_2002_on}, that the total variation distance upper bounds the Lévy–Prokhorov distance \cite{gibbs_2002_on}, the fact that the total variation distance between laws of random variables is invariant to shifting and rescaling the random variables, equation (\ref{TV vertex urn}) and Theorem \ref{Urn Gauss} will complete the proof in the continuous time case.

\begin{align*}
    &\pi\left(\mathcal{L}\left[\frac{|V_t|-\mathbb{E}[|V_t|]}{\sqrt{\mathbb{V}[|U_t|]}}\right], \mathcal{N}(0,1)\right) \\ 
    &\leq \pi\left(\mathcal{L}\left[\frac{|V_t|-\mathbb{E}[|V_t|]}{\sqrt{\mathbb{V}[|U_t|]}}\right], \mathcal{L}\left[\frac{|U_t|-\mathbb{E}[|U_t|]}{\sqrt{\mathbb{V}[|U_t|]}}\right] \right) +\pi\left(\mathcal{L}\left[\frac{|U_t|-\mathbb{E}[|U_t|]}{\sqrt{\mathbb{V}[|U_t|]}}\right], \mathcal{N}(0,1)\right) \\ 
    &\leq d_{TV}\left(\mathcal{L}\left[\frac{|V_t|-\mathbb{E}[|V_t|]}{\sqrt{\mathbb{V}[|U_t|]}}\right], \mathcal{L}\left[\frac{|U_t|-\mathbb{E}[|U_t|]}{\sqrt{\mathbb{V}[|U_t|]}}\right] \right) +\pi\left(\mathcal{L}\left[\frac{|U_t|-\mathbb{E}[|U_t|]}{\sqrt{\mathbb{V}[|U_t|]}}\right], \mathcal{N}(0,1)\right) \\
    &= d_{TV}\left(\mathcal{L}[|V_t|],\mathcal{L}[|U_t|]\right)  +\pi\left(\mathcal{L}\left[\frac{|U_t|-\mathbb{E}[|U_t|]}{\sqrt{\mathbb{V}[|U_t|]}}\right], \mathcal{N}(0,1)\right) \rightarrow 0.
    \end{align*}

To de-Poissonize one uses essentially the same argument as in Theorem 2 of \cite{dutko_1989_central}. We provide the argument in the appendix.
\end{proof}
\begin{remark}\label{connected not very restrictive}
    We interpret Example \ref{gamma > 2} as suggesting that the assumption of eventual forever connectedness is not very restrictive.
\end{remark}

Certainly eventual forever connectedness is not necessary for asymptotic Gaussianity of the vertex count, as the following proposition illustrates.

\begin{proposition} \label{sparse Gaussian} Let the support of $\mu$ consist only of isolated edges and assume that $\mathbb{V}[|V_t|] \rightarrow \infty$ or equivalently $\mathbb{V}[|E_t|] \rightarrow \infty$. Then
$$\frac{|V_t|-\mathbb{E}[|V_t|]}{\sqrt{\mathbb{V}[|V_t|]}} \rightarrow \mathcal{N}(0,1).$$
\end{proposition}
\begin{proof}
Simply noting that the number of vertices is just twice the number of edges gives
$$\frac{|V_t|-\mathbb{E}[|V_t|]}{\sqrt{\mathbb{V}[|V_t|]}} = \frac{2|E_t|-2\mathbb{E}[|E_t|]}{\sqrt{4\mathbb{V}[|E_t|]}} \rightarrow \mathcal{N}(0,1),$$ by Theorem \ref{Urn Gauss}.
\end{proof}

In the proof of Theorem \ref{connected Gauss} the scaling that naturally emerges is that involving the urn scheme. This variance is convenient because it is straightforward to compute. However, one may desire a theorem statement which is in terms of only the vertex process. We will allow this by showing that in the connected regime the two variances are within an additive constant and hence asymptotically equivalent.

\begin{proposition}\label{equi var ineq}
$$\mathbb{V}[|U_t|] \leq \mathbb{V}[|V_t|] \leq \mathbb{V}[|U_t|]+2\sum_{e}\frac{\mu_e}{M_e}.$$ 
\end{proposition}
\begin{proof}  Let $X_i$ be the event that $i \in V_t$ so that $|V_t|=\sum_i 1_{X_i}$. 
    \begin{align*}
        \mathbb{V}[|V_t|] &= \mathbb{V}\left[\sum_{i}1_{X_i}\right], \\
        &= \sum_i \mathbb{V}[1_{X_i}] + \sum_{i \neq j}P(X_i \wedge X_j)-P(X_i)P(X_j), \\
        &= \mathbb{V}[|U_t|]+\sum_{i \neq j}P(X_i \wedge X_j)-P(X_i)P(X_j),
    \end{align*} 
    where in the last equality we recalled that the urns are independent and that $P(i \in V_t) = 1-e^{-M_it}= P(i \in U_t)$.
For compactness of notation, let $M_{i \backslash j} := M_{i}-\mu_{ij}$.
\begin{align*}
    P(X_i \wedge X_j) &= P(X_i \wedge X_j|\tau_{ij} \leq t)P(\tau_{ij} \leq t)+P(X_i \wedge X_j | \tau_{ij} > t)P(\tau_{ij} > t), \\
    &= 1-e^{-\mu_{ij}t}+(1-e^{-M_{i \backslash j}t})(1-e^{-M_{j\backslash i}t})e^{-\mu_{ij}t}, \\
    &= 1-e^{-M_it}-e^{-M_jt}+e^{-M_{ij}t}.
\end{align*}
Since $$P(X_i)P(X_j) = (1-e^{-M_it})(1-e^{-M_jt}),$$ straightforward algebra gives \begin{align*}
    P(X_i \wedge X_j)-P(X_i)P(X_j) = e^{-M_{ij}t}(1-e^{-\mu_{ij}t}).
\end{align*} We note that this is non-negative. Some elementary calculus gives $e^{-at}(1-e^{-bt}) \leq \frac{b}{a}$ for positive reals $a,b,t$ and applying this gives
$$0 \leq P(X_i \wedge X_j)-P(X_i)P(X_j) \leq \frac{\mu_{ij}}{M_{ij}}.$$ Now summation, with the observation that every edge contributes twice to the sum, completes the proof.
\end{proof}

\begin{corollary}
    If $\mu$ is such that $G_t$ is eventually forever connected and at least one of $\mathbb{V}[|U_t|]$ and $\mathbb{V}[|V_t|]$ goes to infinity, then so does the other and \begin{align*}
        \frac{\mathbb{V}[|V_t|]}{\mathbb{V}[|U_t|]} \rightarrow 1.
    \end{align*}
\end{corollary}
\begin{proof}
    By Theorem \ref{Eventual forever connectedness thm} eventual forever connectedness implies that $\sum_e \frac{\mu_e}{M_e} < \infty$. The result then follows from Proposition \ref{equi var ineq} by dividing by $\mathbb{V}[|U_t|]$ and taking limits.
\end{proof}

\begin{remark} 
    The proof of Proposition \ref{equi var ineq} gives some more intuition for why one might expect the connectedness assumption to help with showing Gaussianity. Namely, the calculations in the proof show that $\frac{\mu_{ij}}{M_{ij}}$, the quantity we have already seen is $P(I_{\{i,j\}})$, also serves as an upper bound on the amount of dependence, as measured by the covariance, between the indicators of the vertices $i$ and $j$.
\end{remark}

\section{Completeness}\label{Completeness}

If the set of vertices with positive marginal intensity is infinite it is not possible for the graphs to be eventually forever complete. This is because when a new vertex is added the graph is necessarily not complete and new vertices are added infinitely often. When the support of $\mu$ is infinite the closest one could hope to get to eventual forever completeness is captured in Definition \ref{essentially complete definition}.

\begin{definition}\label{essentially complete definition} We say that a graph $G=(V,E)$ is \textbf{essentially complete} if it is connected and there is an $n \in \mathbb{N}$ such that $V=\{1,\ldots,n+1\}$ and the induced subgraph obtained by taking the vertex subset $\{1,\ldots,n\}$ is the complete graph on $\{1,\ldots,n\}$. 
\end{definition}
We illustrate this definition by giving some examples and non-examples in Figure \ref{fig:essential completesness}.
\begin{figure}[t]
    \centering

\begin{tikzpicture}[every node/.style={draw, circle, fill=blue!30}]
    \node (1a) at (0, 0) {1};
    \node (2a) at (1.5, 0) {2};
    \node (3a) at (0.75, 1.3) {3};
    \node (4a) at (0.75, -1.3) {4};
    \node (5a) at (0.75, 2.6) {5};
    
    \draw (1a) -- (2a);
    \draw (1a) -- (3a);
    \draw (1a) -- (4a);
    \draw (1a) -- (5a);
    \draw (2a) -- (4a);
    \draw (2a) -- (5a);
    \draw (3a) -- (4a);
    \draw (3a) -- (5a);

    \node (1b) at (3.5, 0) {1};
    \node (2b) at (5, 0) {2};
    \node (3b) at (4.25, 1.3) {3};
    \node (4b) at (4.25, -1.3) {4};
    \node (5b) at (4.25, 2.6) {6};
    
    \draw (1b) -- (2b);
    \draw (1b) -- (3b);
    \draw (1b) -- (4b);
    \draw (1b) -- (5b);
    \draw (2b) -- (3b);
    \draw (2b) -- (4b);
    \draw (2b) -- (5b);
    \draw (3b) -- (4b);
    \draw (3b) -- (5b);
    \draw (4b) -- (5b);

    \node (1c) at (7, 0) {1};
    \node (2c) at (8.5, 0) {2};
    \node (3c) at (7.75, 1.3) {3};
    \node (4c) at (7.75, -1.3) {4};
    \node (5c) at (7.75, 2.6) {5};
    
    \draw (1c) -- (2c);
    \draw (1c) -- (3c);
    \draw (1c) -- (4c);
    \draw (2c) -- (3c);
    \draw (2c) -- (4c);
    \draw (2c) -- (5c);
    \draw (3c) -- (4c);
\end{tikzpicture}
    \caption{Examples and non-examples of essentially complete graphs. The first two are not essentially complete as they fail the induced subgraph condition. The former is missing edges, the latter a vertex (and its associated edges). The third is an essentially complete graph because the subgraph induced by vertices $1,2,3,4$ is the complete graph on four vertices, the fifth vertex is the only extra vertex, and that vertex is not isolated.}
    \label{fig:essential completesness}
\end{figure}
\begin{remark}
    Clearly eventual forever essential completeness is a (much) stronger notion than eventual forever connectedness.
\end{remark}
In our terminology, with loops allowed, Janson \cite{janson_2017_on} showed the following.
\begin{proposition}(\cite{janson_2017_on}, Example 8.1) \label{complete sufficient} If for all $k \geq 2$,
$$0 < \sup_{l}\mu_{\{k+1,l\}} \leq k^{-4}\min_{i<k}\mu_{\{k,i\}},$$ then $G_t$ is eventually forever essentially complete. Moreover, $\mu_{ij} \propto (\max(i,j)!)^{-4}$ and $\mu_{ij} \propto e^{-3^i}e^{-3^j}$ satisfy this assumption.
\end{proposition}
We build up to giving a characterization of eventual forever essential completeness.
\begin{definition}
Let $\mathcal{F} =\{\tau_i\}_{i \in \mathcal{I}}$ be a family of random variables. Let $\mathcal{C} =\{C_n\}_{n \in \mathbb{N}}$ be an ordered partition of $\mathcal{I}$ with each $C_n$ nonempty. We define
\begin{align*}
    R := \cup_{n = 1}^\infty \cap_{N=n}^\infty \{\check{C}_N < \hat{C}_{N+1}\},
\end{align*} to be the event that $\mathcal{F}$ \textbf{respects} $\mathcal{C}$.
\end{definition}

\begin{lemma}\label{respect} With notation as above,
    let $\mathcal{F}$ consist of independent exponential random variables with positive intensities $\{\lambda_i\}_{i \in \mathcal{I}}$. Let $\Lambda_ n := \sum_{i \in C_n} \lambda_i$. Then $\mathcal{F}$ respects $\mathcal{C}$ with probability one if
    \begin{align*}
        \sum_n 1-\int_0^\infty \prod_{i \in C_n}(1-e^{-\lambda_it}) \Lambda_{n+1}e^{-\Lambda_{n+1}t}dt < \infty.
    \end{align*} Else $\mathcal{F}$ respects $\mathcal{C}$ with probability zero.
\end{lemma}
\begin{proof}
    Let $D_n := \{\check{C}_n < \hat{C}_{n+1}\}$. Then $R$ is precisely the event that only finitely many $\overline{D}_n$ occur. Since the CDF of a maximum of independent random variables is the product of the CDFs and the minimum of independent exponential random variables is again exponential with intensity being the sum of the intensities and because $\check{C}_n$ and $\hat{C}_{n+1}$ are measurable with respect to independent sigma algebras,
    \begin{align*}
        P(D_n) =  \int_0^\infty \prod_{i \in C_n}(1-e^{-\lambda_it}) \Lambda_{n+1}e^{-\Lambda_{n+1}t}dt.
    \end{align*}
    By the first Borel-Cantelli lemma, if \begin{align*}
        \sum_n 1-P(D_n) = \sum_n P(\overline{D}_n) < \infty,
    \end{align*}
    then $R$ happens almost surely. Suppose now that $\sum_n P(\overline{D}_n) = \infty$. Then either $\sum_{n}P(\overline{D}_{2n}) = \infty$ or $\sum_n P(\overline{D}_{2n+1}) = \infty$. Notice that $\overline{D}_{2n}$ and $\overline{D}_{2m}$ are independent for $n \neq m$. This is because $D_{2n}$ is measurable with respect to the sigma algebra generated by the random variables in $C_{2n} \cup C_{2n+1}$ while $D_{2m}$ is measurable with respect to the random variables in $C_{2m} \cup C_{2m+1}$ and these sigma algebras are independent. Similarly, $D_{2n+1}$ and $D_{2m+1}$ are independent for $n \neq m$. Thus if $\sum_n P(\overline{D}_n) = \infty$ we have a collection of independent events with probabilities that sum to $\infty$. By the second Borel-Cantelli lemma almost surely infinitely many of them will occur. This implies $P(R) = 0$.
\end{proof}

The following result regarding the urn schemes is a simple corollary of Lemma \ref{respect}.

\begin{corollary} In the infinite urn scheme with positive intensities or probabilities, $\{\lambda_i\}_{i \in \mathbb{N}}$, urns are eventually forever filled in order if and only if,
\begin{align*}
    \sum_n \frac{\lambda_{n+1}}{\lambda_n} < \infty.
\end{align*}
\end{corollary}
\begin{proof} By an argument similar to Lemma \ref{De-Poissonization Lemma} it suffices to consider the continuous time urn scheme. With partitions being singletons we get
    \begin{align*}
        P(\overline{D}_n) &= 1-\int_0^\infty (1-e^{-\lambda_n t})\lambda_{n+1}e^{-\lambda_{n+1}t}dt, \\
        &=\frac{\lambda_{n+1}}{\lambda_n+ \lambda_{n+1}}.
    \end{align*}
By the limit comparison test and a little bit of extra work one sees that $\sum_n  \frac{\lambda_{n+1}}{\lambda_n+\lambda_{n+1}}$  is equiconvergent with $\sum_n \frac{\lambda_{n+1}}{\lambda_n}$. 
\end{proof}
\begin{example}
    With $\lambda_n = \frac{1}{n!^\gamma}$ for some $\gamma > 0$ the urns are eventually forever filled in order if and only if $\gamma > 1$. This is clear since $\frac{\lambda_{n+1}}{\lambda_n} = \frac{1}{(n+1)^\gamma}$.
\end{example}

We turn to our main motivation behind Lemma \ref{respect}. Namely, to answer Janson's question about completeness.

\begin{theorem}

Let $\Lambda_{n} = \sum_{i,j: \max(i,j) = n}\mu_{ij}$. If the support of $\mu$ is infinite, then $G_t$ and $G_n$ are eventually forever essentially complete if and only if the support of $\mu$ is the complete graph on $\mathbb{N}$ and 
    \begin{align*}
        \sum_n 1-\int_{0}^\infty \prod_{i,j: \max(i,j)=n}(1-e^{-\mu_{ij}t})\Lambda_{n+1} e^{-\Lambda_{n+1}t}dt < \infty.
    \end{align*}If the support of $\mu$ is finite then instead the support being the complete graph on some vertex set is both necessary and sufficient for $G_t$ and $G_n$ to be eventually forever complete.
\end{theorem}
\begin{proof}
As usual it suffices to consider $G_t$. The finite support case is trivial, so we may focus on $\mu$ with infinite support. If the support is not the complete graph on $\mathbb{N}$, pick an edge missing from the support. Wait until a vertex is added to the graph which is larger, in the sense of the usual ordering of the naturals, than either endpoint of the missing edge. The absence of this edge rules out essential completeness beyond this point. Thus we may restrict to the case that the support is the complete graph on $\mathbb{N}$.

We will now show that eventual forever essential completeness is, up to an event with probability zero, equivalent to the edge arrival times respecting the partition $\{i,j\} \in C_n \Leftrightarrow  \max(i,j) = n$.

Suppose first that the partition is respected. That is, there is some $n_0$ so that $\check{C}_n < \hat{C}_{n+1}$ for all $n \geq n_0$. Since $\mu_e > 0$ and there are only finitely many $\tau_e$ in each block, eventually all edges in blocks that precede $C_{n_0}$ will have appeared. Call this time $t_0$. After $t_0$ there is a number of blocks from which all edges have arrived, then at most one block with some but not all edges having arrived, and then blocks in which no edges have arrived. The first class corresponds to the complete induced subgraph. The second corresponds to the potential extra connected vertex.

Suppose now that the partition is not respected. Then there exist infinitely many $n$ with $\check{C}_{n} > \hat{C}_{n+1}$. For such an $n$, at time $\hat{C}_{n+1}$ the vertex set includes $n+1$, because an edge just brought it. On the other hand, there is an edge between $n$ and some $m < n$ which is missing, so the completeness of the induced subgraph cannot be satisfied. The sequence of such $\hat{C}_{n+1}$ times must be unbounded, since otherwise infinitely many edges would have to arrive in finite time.

Having established the equivalence, and recalling that we may still assume that the support is the complete graph, applying Lemma \ref{respect} gives the result.
\end{proof}

\begin{example} \label{completeness example}
    Take $\mu_{ij} \propto (\max(i,j)!)^{-\gamma}$ for some $\gamma > 0$. Then $G_t$ and $G_n$ are eventually forever essentially complete if and only if $\gamma > 2$.
\end{example}    
    
\begin{proof}
    As before we may assume that the proportionality is in fact an equality. For compactness of notation, let $q_n = (n!)^{-\gamma}$. Then, with more details in the appendix,
    \begin{align*}
        &\sum_n 1-\int_{0}^\infty \prod_{i,j: \max(i,j)=n}(1-e^{-\mu_{ij}t})\Lambda_{n+1} e^{-\Lambda_{n+1}t}dt, \\
        &= \sum_n 1-\int_0^\infty (1-e^{-q_nt})^{n-1}nq_{n+1}e^{-nq_{n+1} t} dt, \\
        &=  \sum_n 1-\prod_{k=1}^{n-1} \frac{k}{\frac{nq_{n+1}}{q_n}+k}, \\
        &= \sum_n1-\prod_{k=1}^{n-1}\frac{k}{\frac{n}{(n+1)^\gamma}+k},
    \end{align*}
    which converges when $\gamma > 2$.
\end{proof}
Example \ref{completeness example} demonstrates that the sufficient condition of \cite{janson_2017_on}, which requires $\gamma \geq 4$ is meaningfully non-sharp.

\printbibliography

\appendix
\section{De-Poissonization}
This is the de-Poissonization required for Theorem \ref{connected Gauss}. It is essentially exactly the same de-Poissonization argument as in \cite{dutko_1989_central}. Let
\begin{align*}
    m(n) &:= \sum_i 1-e^{-M_i n}, \\
    m_n &:= \sum_i 1-(1-M_i)^n, \\
    \sigma^2(t) &:= \mathbb{V}[|U_t|] = m(2t)-m(t).
\end{align*}

Throughout we will use the standing assumption that $\sigma(t) \rightarrow \infty$, implying in particular that there are infinitely many non-zero $M_i$.

\begin{lemma} \label{lem: D-1}
    \begin{align*}
        m(n)-m_n \rightarrow 0.
    \end{align*}
\end{lemma}
\begin{proof} As in \cite{dutko_1989_central} we use that for $0 \leq x \leq n$ ,
    \begin{align*}
        0 \leq e^{-x}-\left(1-\frac{x}{n}\right)^n \leq \frac{x^2}{n}e^{-x}.
    \end{align*}
We set $x = nM_i \leq n$ and sum over $i$ to get
\begin{align*}
\sum_i e^{-M_in}-(1-M_i)^n \leq \sum_i nM_i^2e^{-M_in} \leq 
\sum_i \frac{1}{e} M_i, \end{align*}
where we used that $ye^{-y} \leq e^{-1}$. The terms in the middle sum all go to zero and since $\sum_i M_i =2 < \infty$ the second inequality supplies a dominating map and so the lemma follows by dominated convergence.
\end{proof}

\begin{lemma} \label{lem: D-4}
    \begin{align*}
        \frac{m'(t)^2}{m'(2t)} \rightarrow 0.
    \end{align*}
\end{lemma}
\begin{proof}
    Pick a finite subset $S \subset \mathbb{N}$. For each $k$ in $S$, find some $l$ so that $M_l < M_k$. This is possible since there are infinitely many positive $M_i$. Notice that
    \begin{align*}
        \frac{M_ke^{-M_kt}}{m'(t)} \leq \frac{M_ke^{-M_kt}}{M_le^{-M_lt}} \rightarrow 0.
    \end{align*}
    Thus the terms in $S$ are asymptotically negligible. Using this, writing $M_ie^{-M_it} = \sqrt{M_i}\sqrt{M_i}e^{-M_it}$ then applying the Cauchy-Schwarz inequality gives,
    \begin{align*}
        \limsup_{t \rightarrow \infty} \frac{m'(t)^2}{m'(2t)} = \limsup_{t \rightarrow \infty} \frac{(\sum_{i \notin S}M_ie^{-M_it})^2}{m'(2t)} \leq \sum_{i \notin S} M_i,
    \end{align*} which can be made arbitrarily small.
\end{proof}
\begin{lemma}\label{lem: D-3}
    For every $C>0$,
\begin{align*}
    \lim_{n \rightarrow \infty} \frac{m(n+C\sqrt{n})-m(n)}{\sigma(n)} = 0.
\end{align*}
\end{lemma}
\begin{proof} The observations $m'(t) = \sum_i M_i e^{-M_it} > 0$ and $m''(t) = -\sum_i M_i^2 e^{-M_it} < 0$ give the bounds
    \begin{align*}
        m(n+C\sqrt{n})-m(n) &= \int_n^{n+C\sqrt{n}}m'(t)dt \leq C \sqrt{n} m'(n), \\
        \sigma^2(n) &= \int_n^{2n}m'(t)dt \geq nm'(2n).
    \end{align*}
From these we have
\begin{align*}
    \frac{m(n+C\sqrt{n})-m(n)}{\sigma(n)} \leq \frac{C \sqrt{n}m'(n)}{\sqrt{n}\sqrt{m'(2n)}} \rightarrow 0,
\end{align*} by Lemma \ref{lem: D-4}.
\end{proof}
\begin{lemma} \label{lem: D-5}
    For every $C > 0$,
    \begin{align*}
        \frac{m_{n+C\sqrt{n}}-m_n}{\sigma(n)} \rightarrow 0.
    \end{align*}
\end{lemma}
\begin{proof} Writing
    \begin{align*}
        \frac{m_{n+C\sqrt{n}}-m_n}{\sigma(n)}=&\frac{m_{n+C\sqrt{n}}-m(n+C\sqrt{n})}{\sigma(n)} \\
        +&\frac{m(n+C\sqrt{n})-m(n)}{\sigma(n)}+\frac{m(n)-m_n}{\sigma(n)},
    \end{align*} and applying Lemma \ref{lem: D-1} to the first and third terms and Lemma \ref{lem: D-3} to the middle term.
\end{proof}
\begin{lemma} \label{lem: D-6}
    For every $C,\epsilon > 0$,
    \begin{align*}
        P(|V_{n+C\sqrt{n}}|-|V_n| > \epsilon \sigma(n)) \rightarrow 0.
    \end{align*}
\end{lemma}
\begin{proof}By Markov's inequality and Lemma \ref{lem: D-5}
\begin{align*}
P(|V_{n+C\sqrt{n}}|-|V_n| > \epsilon \sigma(n)) \leq \frac{m_{n+C\sqrt{n}}-m_n}{\epsilon \sigma(n)} \rightarrow 0.
\end{align*}
\end{proof}

Now we are ready to carry out the de-Poissonization in Theorem \ref{connected Gauss}.
\begin{proof}
Let
    \begin{align*}
        \Phi_n(x) := P\left(\frac{|V_t|-m_n}{\sigma(n)} \leq x\right)\Big|_{t=n}
    \end{align*}
From the continuous time CLT, Slutsky's Theorem and Lemma \ref{lem: D-1}, $\Phi_n(x) \rightarrow \Phi(x)$.

 Let $F_n(x):= P(|V_n|\leq m_n+x\sigma(n))$ and $p_k(n):= \frac{n^k}{k!}e^{-n}$. Fix $\delta,\epsilon > 0$. Choose $C$ such that 
 \begin{align*}
     \sum_{k:|k-n| > C \sqrt{n}}p_k(n) < \epsilon,
 \end{align*}
for $n$ large enough. For such large $n$, consider $k$ such that $|k-n| \leq C\sqrt{n}$.

$k \geq n$, $|V_k| \geq |V_n|$ gives $P(|V_k| \leq m_n+x\sigma(n)) \leq F_n(x)$. To find a nearly matching lower bound, note that on the event,
\begin{align*}
    \{|V_n| \leq m_n + (x-\delta)\sigma(n)\} \cap \{|V_k|-|V_n| \leq \delta \sigma(n)\},
\end{align*}
we have $|V_k| \leq m_n + x\sigma(n)$ implying, by Lemma \ref{lem: D-6}, that for $n$ large enough, $$P(|V_k| \leq m_n+x\sigma(n)) \geq F_n(x-\delta)-\epsilon.$$
$k$ close to but smaller than $n$ are handled similarly. This gives that for all $k$ such that $|k-n| \leq C \sqrt{n}$,
\begin{align*}
    F_n(x-\delta)-\epsilon \leq P(|V_k| \leq m_n+x\sigma(n)) \leq F_n(x+\delta)+\epsilon.  
\end{align*}

\begin{align*}
    \Phi_n(x) &= \sum_{k}P(|V_k| \leq m_n+x\sigma(n))p_k(n) \\
    &\leq \sum_{k: |k-n| > C\sqrt{n}}p_k(n)+\sum_{k: |k-n| \leq C\sqrt{n}}(F_n(x+\delta)+\epsilon)p_k(n) \\
    &\leq F_n(x+\delta)+2\epsilon.
\end{align*}
Using the lower bound similarly gives
\begin{align*}
    F_n(x-\delta)-2\epsilon \leq \Phi_n(x) \leq F_n(x+\delta)+2\epsilon.
\end{align*}
Recalling that we already established that $\Phi_n(x) \rightarrow \Phi(x),$ where the latter is the cdf of the standard normal, this implies

\begin{align*}
    \Phi(x-\delta)\leq \liminf_n F_n(x) \leq \limsup_n F_n(x) \leq \Phi(x+\delta),
\end{align*}
which gives the theorem by sending $\delta$ to $0$ and using continuity of $\Phi(x)$.
\end{proof}

\section{Asymptotics for Completeness Example}
We make the substitution $u = e^{-q_nt}$. Let $r_n := \frac{nq_{n+1}}{q_{n}}$.
\begin{align*}
    I_n :=& \int_0^\infty (1-e^{-q_nt})^{n-1}q_{n+1}ne^{-nq_{n+1}t} dt,\\
    =& r_n\int_{0}^1 (1-u)^{n-1}u^{r_n-1}du, \\
    =& r_n B(r_n,n), \\
\end{align*} where $B$ denotes the Beta-function. Now writing $B(r_n,n) = \frac{\Gamma(r_n)\Gamma(n)}{\Gamma(r_n+n)}$ and repeatedly applying the recursive identity for the gamma function,

\begin{align*}
    I_n = \frac{\Gamma(r_n+1)\Gamma(n)}{\Gamma(r_n+n)} = \prod_{k=1}^{n-1}\frac{k}{r_n+k}.
\end{align*}
Taking logarithms and applying the elementary inequality $x - \frac{x^2}{2} \leq \log(1+x) \leq x$ gives
\begin{align*}
\sum_{k=1}^{n-1} \frac{r_n}{k}-\left(\frac{r_n}{\sqrt{2}k}\right)^2\leq -\log(I_n) \leq \sum_{k=1}^{n-1} \frac{r_n}{k}
\end{align*}
In our particular case $r_n = \frac{n(n!)^\gamma}{(n+1)!^\gamma} = \frac{n}{(n+1)^\gamma},$ 
so for $\gamma > 1$, $r_n \rightarrow 0$ (For $\gamma \leq 1, I_n \nrightarrow 1$ and there is no hope of summability) and so for $n$ large enough we have

\begin{align*}
\frac{1}{2}H_{n-1}r_n \leq -\log(I_n) \leq H_{n-1}r_n,
\end{align*} where $H_n$ is the $n$th harmonic number. Applying the increasing function $1-e^{-x}$ to both sides and applying the estimate $x-\frac{x^2}{2} \leq 1-e^{-x} \leq x$ gives

\begin{align*}
    \frac{1}{2}r_n H_{n-1}-\frac{1}{8}(r_nH_{n-1})^2 \leq 1- I_n \leq H_{n-1}r_n.
\end{align*}
Since the harmonic numbers grow logarithmically in $n$, for $\gamma>1$ we have $r_n H_{n-1} \rightarrow 0$. Thus for $n$ large enough we have 
\begin{align*}
    \frac{1}{4}r_n H_{n-1}\leq 1- I_n \leq r_nH_{n-1}.
\end{align*}
Again because the harmonic numbers grow logarithmically it suffices to consider $\sum_n \log(n)n^{1-\gamma}$, which converges precisely when $\gamma > 2$.

\end{document}